\documentclass[12pt,notitlepage,reqno]{amsart} 
 \usepackage{amsmath}     
 \makeatletter 
  \def\activeat#1{\csname @#1\endcsname} 
  \def\def@#1{\expandafter\def\csname @#1\endcsname} 
  {\catcode`\@=\active \gdef@{\activeat}} 

 \let\ssize\scriptstyle 
 \newdimen\ex@        \ex@.2326ex 

  \def\requalfill{\cleaders\hbox{$\mkern-2mu\mathord=\mkern-2mu$}\hfill 
   \mkern-6mu\mathord=$} 
  \def\eqfill{$\m@th\mathord=\mkern-6mu\requalfill} 
  \def\deffill{\hbox{$:=$}$\m@th\mkern-6mu\requalfill} 
  \def\fiberbox{\hbox{$\vcenter{\hrule\hbox{\vrule\kern1ex 
      \vbox{\kern1.2ex}\vrule}\hrule}$}} 

  \newdimen\arrwd 
  \newdimen\minCDarrwd \minCDarrwd=2.5pc 
           
  \def\findarrwd#1#2#3{\arrwd=#3%
   \setbox\z@\hbox{$\ssize\;{#1}\;\;$}%
   \setbox\@ne\hbox{$\ssize\;{#2}\;\;$}%
   \ifdim\wd\z@>\arrwd \arrwd=\wd\z@\fi 
   \ifdim\wd\@ne>\arrwd \arrwd=\wd\@ne\fi} 
  \newdimen\arrowsp\arrowsp=0.375em           
  \def\findCDarrwd#1#2{\findarrwd{#1}{#2}{\minCDarrwd} 
     \advance\arrwd by 2\arrowsp} 
  \newdimen\minarrwd  
  \setbox\z@\hbox{$\longrightarrow$} \minarrwd=\wd\z@ 

  \def\harrow#1#2#3#4{{\minarrwd=#1\minarrwd 
    \findarrwd{#2}{#3}{\minarrwd}\kern\arrowsp 
     \mathrel{\mathop{\hbox to\arrwd{#4}}\limits^{#2}_{#3}}\kern\arrowsp}} 

  \def@]#1>#2>#3>{\harrow{#1}{#2}{#3}\rightarrowfill} 
  \def@>#1>#2>{\harrow1{#1}{#2}\rightarrowfill} 
  \def@<#1<#2<{\harrow1{#1}{#2}\leftarrowfill} 
  \def@={\harrow1{}{}\eqfill} 
  \def@:#1={\harrow1{}{}\deffill} 
  \def@ N#1N#2N{\vCDarrow{#1}{#2}\UpDownarrow} 
  \def\UpDownarrow{\uparrow\,\Big\downarrow} 

 \def\hookrightarrowfill{\hbox{$\lhook\joinrel$}\rightarrowfill} 
 \def@{hk>}#1>#2>{\harrow1{#1}{#2}\hookrightarrowfill} 
 \def@ c>#1>#2>{\harrow1{#1}{#2}\hookrightarrowfill} 
 \def\hookleftarrowfill{\leftarrowfill\hbox{$\joinrel\rhook$}} 
 \def@{hk<}#1<#2<{\harrow1{#1}{#2}\hookleftarrowfill} 

  \def@.{\ifodd\row\relax\harrow1{}{}\hfill 
    \else\vCDarrow{}{}.\fi} 
  \def@|{\vCDarrow{}{}\Vert} 
  \def@ V#1V#2V{\vCDarrow{#1}{#2}\downarrow} 
  \def@ A#1A#2A{\vCDarrow{#1}{#2}\uparrow} 
  \def@(#1){\arrwd=\csname col\the\col\endcsname\relax 
    \hbox to 0pt{\hbox to \arrwd{\hss$\vcenter{\hbox{$#1$}}$\hss}\hss}} 

  \def\squash#1{\setbox\z@=\hbox{$#1$}\finsm@@sh} 
 \def\finsm@@sh{\ifnum\row>1\ht\z@\z@\fi \dp\z@\z@ \box\z@} 

  \newcount\row \newcount\col \newcount\numcol \newcount\arrspan 
  \newdimen\vrtxhalfwd  \newbox\tempbox 

  \def\innernewdimen{\alloc@1\dimen\dimendef\insc@unt} 
  \def\measureinit{\col=1\vrtxhalfwd=0pt\arrspan=1\arrwd=0pt 
    \setbox\tempbox=\hbox\bgroup$} 
  \def\setinit{\col=1\hbox\bgroup$\ifodd\row 
    \kern\csname col1\endcsname 
    \kern-\csname row\the\row col1\endcsname\fi} 
  \def\findvrtxhalfsum{$\egroup 
   \expandafter\innernewdimen\csname row\the\row col\the\col\endcsname 
   \global\csname row\the\row col\the\col\endcsname=\vrtxhalfwd 
   \vrtxhalfwd=0.5\wd\tempbox 
   \global\advance\csname row\the\row col\the\col\endcsname by \vrtxhalfwd 
   \advance\arrwd by \csname row\the\row col\the\col\endcsname 
   \divide\arrwd by \arrspan 
   \loop\ifnum\col>\numcol \numcol=\col%
      \expandafter\innernewdimen \csname col\the\col\endcsname 
      \global\csname col\the\col\endcsname=\arrwd 
    \else \ifdim\arrwd >\csname col\the\col\endcsname 
       \global\csname col\the\col\endcsname=\arrwd\fi\fi 
    \advance\arrspan by -1 %
    \ifnum\arrspan>0 \repeat} 
  \def\setCDarrow#1#2#3#4{\advance\col by 1 \arrspan=#1 
     \arrwd= -\csname row\the\row col\the\col\endcsname\relax 
     \loop\advance\arrwd by \csname col\the\col\endcsname 
      \ifnum\arrspan>1 \advance\col by 1 \advance\arrspan by -1%
      \repeat 
     \squash{\mathop{ 
      \hbox to\arrwd{\kern\arrowsp#4\kern\arrowsp}}\limits^{#2}_{#3}}} 
  \def\measureCDarrow#1#2#3#4{\findvrtxhalfsum\advance\col by 1%
    \arrspan=#1\findCDarrwd{#2}{#3}%
     \setbox\tempbox=\hbox\bgroup$} 
  \def\vCDarrow#1#2#3{\kern\csname col\the\col\endcsname 
     \hbox to 0pt{\hss$\vcenter{\llap{$\ssize#1$}}%
      \Big#3\vcenter{\rlap{$\ssize#2$}}$\hss}\advance\col by 1} 

  \def\setCD{\def\harrow{\setCDarrow}%
   \def\\{$\egroup\advance\row by 1\setinit} 
   \m@th\lineskip3\ex@\lineskiplimit3\ex@ \row=1\setinit} 
  \def\endsetCD{$\egroup} 
  \def\drop#1\\{\findvrtxhalfsum\advance\row by 2 \measureinit} 
  \def\measure{\bgroup 
   \def\harrow{\measureCDarrow}%
   \def\\##1{\ifx##1\endmeasure\endmeasure\else\expandafter\drop\fi}%
   \row=1\numcol=0\measureinit} 
  \def\endmeasure{\findvrtxhalfsum\egroup} 

  \newcount\savedcount         
  \def\LCD#1\end{\savedcount=\count11 
    \measure#1\endmeasure 
    \vcenter{\setCD#1\endsetCD\kern\medskipamount}%
    \global\count11=\savedcount\end} 

  \newenvironment{CD}{\let\at=@\catcode`\@=\active\LCD}{\catcode`\@=12\relax} 
 \makeatother 

 \usepackage{amsfonts,amsthm} 
 \usepackage{amssymb} 

 \advance\voffset by -0.4 truein 
 \advance\hoffset by -0.75 truein 
 \font\smallrm=cmr8 

 \newcommand{\Og}{\Omega} 
 \renewcommand{\H}{\operatorname{\rm H}} 
 \newcommand{\h}{\operatorname{\rm h}} 
 \newcommand{\Hom}{\operatorname{\rm Hom}} 
 \newcommand{\Ker}{\operatorname{\rm Ker}} 
 \newcommand{\Cok}{\operatorname{\rm Cok}} 
 \renewcommand{\Im}{\operatorname{\rm Im}} 
 \newcommand{\ox}{\otimes} 
 \newcommand{\CIP}{\text{\bf P}_{\text{\bf C}}} 
 \newcommand{\IP}{\text{\bf P}} 
 \newcommand{\X}{\IP^2} 
 \newcommand{\reg}{\mathrm{reg}\,} 
 \newcommand{\ve}{\varepsilon} 
 \renewcommand{\:}{\colon} 

 \newcommand{\smashedrightarrow} 
  {\setbox0=\hbox{$\longrightarrow$}\ht0=1pt\box0} 
 \newcommand{\risom} 
  {\buildrel{\hskip-0.08cm\sim}\over{\smashedrightarrow}} 

 \newcommand{\vf}{\varphi} 

 \def~{\thinspace\nobreak} 

 \newtheorem{theorem}{Theorem}[section] 
 \newtheorem{lemma}[theorem]{Lemma} 
 \newtheorem{proposition}[theorem]{Proposition} 
 \theoremstyle{definition} 
 \newtheorem{remark}[theorem]{Remark} 
 \newtheorem{subsct}[theorem]{} 
 \theoremstyle{plain} 
   
\begin{document} 
 \title[{\smallrm Bounds on leaves of foliations of the plane}] 
 {Bounds on leaves of foliations of the plane} 
 \author[{\smallrm E. Esteves and S. Kleiman}]{E. Esteves$^1$ 
 \ and \ S. Kleiman$^2$} 
  \date{9 January 2003} 
 \begin{abstract} 
 This paper contributes to the solution of the Poincar\'e problem, which
is to bound the degree of a (generalized algebraic) leaf of a (singular
algebraic) foliation of the complex projective plane.  The first theorem
gives a new sort of bound, which involves the Castelnuovo--Mumford
regularity of the singular locus of the leaf.  The second theorem gives
a bound in terms of two singularity numbers of the leaf: the total
Tjurina number, and the number of non-quasi-homogeneous singularities.
If such singularities are present, then this bound improves one due to du
Plessis and Wall, at least when the curve is irreducible.
 \end{abstract}

\maketitle 

 \section{Introduction} 

When does a singular foliation of $\CIP^2$ defined by polynomials have a 
leaf defined by a polynomial?  This question is fundamental, but difficult, 
and it has 
 stimulated a lot of research for well over a century. 
 \renewcommand{\thefootnote}{}%
 \footnote 
  {2000 {\it Mathematics Subject Classification} 37F75 (primary), 14H50, 
    32S65, 14H20 (secondary).} 
 \footnote 
  {{\it Keywords} Foliations, curves, singularities.} 
 \footnote 
  {$^1$This author thanks A. Campillo, L. G. Mendes, P. Sad, 
 M. Soares, and especially J. V. Pereira for helpful discussions on the 
 subject.  
He is also grateful to CNPq for a grant, Proc.\ 202151/90-5, supporting
a year-long visit to MIT, and grateful to MIT for its hospitality. 
He was also supported by PRONEX, Conv\^enio 41/96/0883/00, CNPq, 
Proc.\ 300004/95-8, and FAPERJ, Proc.{} E-26/170.418/2000-APQ1.
} 
 \footnote 
  {$^2$This author thanks IMPA, Rio de Janeiro, ICMC-USP, S\~ao Carlos,
    and the XVI and XVII Escolas de \'Algebra, Bras\'\i lia and Cabo Frio,
    for their invitations and financial support, which enabled this work
    to be initiated, pursued, and presented.}

 Notably, in 1891, Poincar\'e \cite{P}, p.~161, observed that, once we 
possess a bound on the degree $d$ of a polynomial $F$ defining a leaf, we 
can try to find $F$ by making purely algebraic computations.  The problem of 
finding such bounds is now known as the {\it Poincar\'e problem}.
The  substantial current interest in it was stimulated by Cerveau and 
 Lins~Neto \cite{CL} in 1991. 

 The available bounds on the degree $d$ of $F$ depend on the degree $m$
of the foliation and on the number and type of the singularities of the
Zariski closure $C$ of the leaf, the curve defined by $F$.
 For example, Campillo and Carnicer \cite{CC} found bounds depending on 
 the topological type as encoded in the Enriques resolution diagram; 
 their results are improved in \cite{EK1}, Theorem~5.3.  Du Plessis and 
 Wall \cite{dPW} found bounds depending on the analytic type as encoded 
 in the Tjurina number; their results are generalized to curves in any 
 projective space in \cite{EK1}, Section 6, via a second approach, and 
 they are improved below via a third approach.  Additional work is cited 
 in \cite{EK1}. 

 The present note introduces a new sort of bound, one that also takes 
 into account the distribution of the singularities as measured by the 
 Castelnuovo--Mumford regularity $\sigma$ of the singular locus of $C$. 
In a similar vein, \cite{Es} introduced the study of the regularity of a 
variety invariant under a singular vector field of projective space.  
And \cite{EK2} 
 introduced the study of the regularity of the intersection of a solution 
 of a Pfaff system with the singular locus of the system. 

 In the present note, the first main result is Theorem~2.5, which asserts 
 this: if $\sigma\le d-2$, then $d\le m+1$; otherwise, $d\le m+1+\rho$ 
 where $\rho:=\sigma-d-2$; furthermore, $d=m+1+\rho$ if $d\ge 2m+2$ and 
 if the foliation has finitely many singularities.  The proof requires 
 Corollary~4.5 of \cite{EK2}, but only its assertion that $\h^1(\mathcal 
 O_C(m-1))=0$ if $\h^1(\Og^1_C)=1$. 

 The regularity $\sigma$ satisfies the following upper bound, given in 
 Lemma~3.1: 
         $$\sigma\le d-2 + (\tau-u)/(d-1)$$ 
 where $\tau$ is the total Tjurina number, the sum of the local Tjurina 
 numbers of $C$, and where $u$ is the number of non-quasi-homogeneous 
 singularities of $C$. 

 Lemma 3.1 and Theorem 2.5 yield the second main result, 
 Theorem~3.2, which asserts that
         $$(d-1)(d-m-1)+u\le\tau.$$ 
  By incorporating $u$, this bound improves that of du Plessis and Wall 
\cite{dPW}, Theorem 3.2, at least when $C$ is irreducible.

 In singularity theory, upper bounds on $\tau$ are also important. 
 Proposition 6.3 of \cite{EK1} asserts a bound for curves in any 
 projective space, if the foliation has finitely many singularities. 
 In the present case, the bound is this: 
         $$\tau\le(d-1)(d-m-1)+m^2.$$ 
 If $d\le2m$, then the bound can be improved, provided $C$ is irreducible 
and is not the closure of a leaf of a foliation of smaller degree than $m$; 
namely, 
 then 
         $$\textstyle\tau\le(d-1)(d-m-1)+m^2-\binom{2m+2-d}{2}.$$ 

At least when $C$ is irreducible, these two bounds agree with those
proved
 by du Plessis and Wall \cite{dPW}, Theorem 3.2.  Both bounds are
 asserted in
 Proposition~3.3 below, and proved via a rather different and more 
 conceptual approach, similar to that of \cite{EK1} and \cite{EK2}. 

Du Plessis and Wall proved, with consummate skill, their lower and upper
bounds on $\tau$ under different, possibly weaker hypotheses, which are
discussed in Remark~3.4 below.  However, seen from the viewpoint of
foliation theory, our hypotheses, stronger or not, do not appear to be
unreasonable.

 All the proofs below are purely algebraic, and work over any 
 algebraically closed field of characteristic 0, not just $\bf C$.  In 
 fact, Section~2 and Proposition~3.3 are set over an algebraically closed 
field of arbitrary characteristic $p$, but the proof of the Theorem~2.5 
requires $p\nmid d$. 

 \section{Bounds on the degree} 

 \begin{subsct} \emph{Foliations of the plane}.  By definition, a 
\emph{singular}
 \emph{foliation} of $\X$ is a nonzero map $\eta\:\Og^1_{\X}\to\mathcal L$ with 
 $\mathcal L$ invertible.  Its \emph{singular locus} is the subscheme 
 $S\subset\X$ with ideal 
              $$\mathcal I_{S,\,\X}:=\Im(\eta\ox\mathcal L^{-1}).\eqno(2.1.1)$$ 

 The singular locus $S$ is never empty.  Indeed, its degree is given by
Formula (3.3.2) below, and plainly never 0.  (This formula was,
according to Poincar\'e \cite{P}, p.~165, known before 1870 to Darboux.)

The singular foliation $\eta$ defines an actual foliation of $\X-S$.
For convenience, let's call a singular foliation simply a {\it foliation}.

 Note $\dim S\le 1$ since $\eta\neq 0$.  If $\dim S=1$, let $D\subset\X$ 
 be the largest curve (effective divisor) contained in $S$.  Then $\eta$ 
 factors through a foliation $\eta'\:\Og^1_{\X}\to\mathcal L(-D)$, and its 
 singular locus is finite.  Thus questions about foliations of $\X$ can 
 often be reduced to the case of foliations with finite singular locus. 

 Let $P\in\X$.  If $P\notin S$, then $\eta(P)\:\Og^1_{\X}(P)\to\mathcal L(P)$ 
 is a surjection of vector spaces; hence, $\eta(P)$ defines a 
 1-dimensional subspace of the tangent space $T_P\,\X$, so a line 
 $L_P\subset\X$ through $P$.  Enumerate the $P$ on a general line 
 $M\subset\X$ with $P\in S$ or $L_P=M$; the total is called the 
 \emph{degree} of $\eta$ and denoted by $\deg\eta$.  More precisely, 
 $\deg\eta$ is the degree of the degeneracy locus of the induced map 
 $(\eta|M,\beta_M)\:\Og^1_{\X}|M\to\mathcal L|M\bigoplus\Og^1_M$ where 
 $\beta_M$ is the natural map, which is the degeneracy locus of 
 $\wedge^2(\eta|M,\beta_M)\:\Og^2_{\X}|M\to\mathcal L|M\ox\Og^1_M$. Plainly, 
 $\deg\eta\ge0$.

 Recall $\deg\Og^1_{\X}=-3$ and $\deg\Og^1_M=-2$.  It follows that 
 $\deg\eta =\deg\mathcal L+1$, or 
         $$\mathcal L\cong\mathcal O_{\X}(m-1)\text{ where }m:=\deg\eta.\eqno(2.1.2)$$ 

Let $C\subset\X$ be a reduced curve.  Call $C$ a \emph{generalized} 
\emph{leaf} of $\eta$ 
 if no component of $C$ lies in $S$ and if $\eta|C$ factors 
 through the natural surjection $\beta_C\:\Og^1_{\X}|C\to\Og^1_C$. 

Assume $C\cap S$ is finite.  Then $C$ is a generalized leaf if and only if, 
for 
 every simple point $P$ of $C$ off $S$, the line $L_P$ is the tangent to 
$C$ at $P$.  Indeed, if $C$ is a generalized leaf, then the induced map 
 $\Og^1_C\to\mathcal L|C$ is, at $P$, a surjection from one invertible sheaf to 
 another, so an isomorphism; hence, $L_P$ is the tangent line. 
 Conversely, suppose $\eta|C$ does not factor through $\beta_C$.  Set $\mathcal 
 K:=\Ker\beta_C$.  Then $\eta|\mathcal K$ is nonzero.  Hence, since $C$ is 
reduced and $\mathcal L$ is invertible, $\eta|\mathcal K$ is surjective at infinitely 
many points of a component of $C$.  So $\eta|\mathcal K$ is surjective at some 
 simple point $P$ of $C$.  Then $P\notin S$, and $L_P$ is not the tangent 
 to $C$ at $P$. 

A generalized leaf is thus a union of actual leaves of the foliation defined 
by $\eta$ on $\X-S$.  For added convenience, from now on, let's call a
generalized leaf simply a {\it leaf}.
 \end{subsct} 

 \begin{subsct} \emph{Singularities of plane curves.}  Let $C\subset\X$ 
 be a reduced curve.  By definition, the \emph{singular locus} of $C$ is 
 the subscheme $\Sigma\subset C$ whose ideal $\mathcal I_{\Sigma,C}$ is the 
 first Fitting ideal of $\Og^1_C$. 

 Set $d:=\deg C$.  Then the ideal $\mathcal I_{\Sigma,C}$ can be computed using 
 the standard presentation 
         $$\mathcal O_C(-d)\to\Og^1_{\X}|C\to\Og^1_C\to 0.\eqno(2.2.1)$$ 
  It yields a map  $\mathcal O_C(-d)\ox\Og^1_C\to\wedge^2\Og^1_{\X}|C$, 
 whence a map $\Og^1_C\to\mathcal O_C(d-3)$.  The latter is an isomorphism off 
 $\Sigma$, and its image is $\mathcal I_{\Sigma,C}(d-3)$. Hence 
         $$\Og^1_C/\hbox{Torsion}\risom\mathcal I_{\Sigma,C}(d-3).\eqno(2.2.2)$$

 Let $P\in \Sigma$. Take local coordinates $x,\,y$ for $\X$ at $P$, and 
 let $f(x,y)=0$ be a local equation for $C$.  By definition, the 
 \emph{Tjurina number} $\tau_P$ of $C$ at $P$ is the colength in $\mathcal 
 O_{\X,\,P}$ of the ideal generated by $f$ and its partials $f_x,\,f_y$. 
 And the \emph{Milnor number} $\mu_P$ of $C$ at $P$ is the colength in 
 $\mathcal O_{\X,\,P}$ of the ideal generated only by $f_x,\,f_y$.  Plainly, both 
 $\tau_P$ and $\mu_P$ are analytic (formal) invariants. 

 By definition, $P$ is a \emph{quasi-homogeneous} singularity of $C$ if, 
 after a suitable analytic change of variables, $f(x,y)$ becomes a 
weighted-homogeneous polynomial.  Saito \cite{Sa}, 
second ``Satz'' on p.~123, proved that, in 
 characteristic zero, $P$ is quasi-homogeneous if and only if 
 $\tau_P=\mu_P$; the result is also discussed in \cite{D}, see 
 Theorem~7.42, p.~120. 

 {}From (2.2.1) and (2.2.2), it is clear that the ideal $\mathcal I_{\Sigma,\,C}$ 
 is generated at $P$ by $f_x,\,f_y$.  Hence the Tjurina number $\tau_P$ 
 is the length of $\mathcal O_\Sigma$ at $P$.  Since $C$ is reduced, $\Sigma$ 
 is finite.  Set $\tau:=\sum_{P\in\Sigma}\tau_P$, and call $\tau$ the 
 \emph{total Tjurina number} of $C$.  Note $\tau=\deg 
 \Sigma$. 

 By definition, the \emph{polar system} of $C$ is the linear system on 
 $\X$ generated by the three partial derivatives of the homogeneous 
 polynomial $F$ defining $C$ (see \cite{Do}, Section~3).  If the 
 characteristic is $0$, or if it is positive and does not divide $\deg 
 C$, then Euler's formula shows that $F$ belongs to the ideal generated 
 by its partial derivatives.  Therefore, the base locus of the polar 
 system is, scheme-theoretically, just the singular locus $\Sigma$ of 
 $C$. 

 \end{subsct} 

 \begin{proposition} Let $\eta$ be a foliation of $\X$.  Let $S$ be its 
 singular locus, and $m$ its degree.  If $S$ is finite and if $m>0$, then 
 $\reg S=2m$. 
 \end{proposition} 

 \begin{proof} Owing to (2.1.1) and (2.1.2), the Koszul complex on 
 $\eta\ox\mathcal L^{-1}$ yields a short sequence 
         $$0\to\Og^2_{\X}(2-2m)\to\Og^1_{\X}(1-m)\to\mathcal I_{S,\,\X}\to 0. 
         \eqno(2.3.1)$$ 
 It is exact since $S$ is finite and $\X$ is Cohen-Macaulay. 

 Twisting (2.3.1) by $i-1$, and taking cohomology, we obtain the 
 following exact sequence: 
  $$\H^1(\Og^1_{\X}(i-m))\to\H^1(\mathcal I_{S,\,\X}(i-1))
   \to\H^2(\Og^2_{\X}(i+1-2m))\to\H^2(\Og^1_{\X}(i-m)).\eqno(2.3.2)$$
  Now, $\H^1(\Og^1_{\X}(i-m))\neq 0$ if and only if $i=m$, by \cite{De},
Th\'eor\`eme~1.1, p.~40.  By duality, $\H^2(\Og^2_{\X}(i+1-2m))=0$ if
and only if $i\geq 2m$.  By hypothesis, $m>0$, or $2m>m$.  Take $i=2m$
in (2.3.2).  Its exactness now yields $\H^1(\mathcal I_{S,\,\X}(2m-1))=0$.
Thus $\reg S\leq 2m$.

 On the other hand, by \cite{De}, Th\'eor\`eme~1.1, p.~40, again, 
 $\H^2(\Og^1_{\X}(m-1))=0$ as $m\geq 0$.  Now, $\H^2(\Og^2_{\X})\neq 0$. 
 Take $i=2m-1$ in (2.3.2).  Its exactness yields $\H^1(\mathcal 
 I_{S,\,\X}(2m-2))\neq0$.  Thus $\reg S\geq 2m$. 
  \end{proof} 
    
 \begin{lemma} 
   Let $C\subset\X$ be a curve of degree $d$, and $T\subset C$ a finite 
 subscheme.  Then 
 \begin{align*} 
  \h^2(\mathcal I_{T,\,\X}(i))&=0\text{ for }i\ge-2,\tag{2.4.1}\\ 
  \h^1(\mathcal I_{T,\,C}(i))&=\h^1(\mathcal I_{T,\,\X}(i))+\epsilon 
  \text{ where } \epsilon := \begin{cases}1,        &\text{if }i=d-3;\\ 
                                         0,        &\text{if }i>d-3. 
                            \end{cases}\tag{2.4.2} 
  \end{align*} 
  \end{lemma} 

 \begin{proof} For any $i$, consider the two twisted standard exact sequences: 
  \begin{gather*}0\to\mathcal I_{T,\,\X}(i)\to\mathcal O_{\X}(i)\to\mathcal O_T(i)\to 0,\\ 
 0\to\mathcal O_{\X}(i-d)\to\mathcal I_{T,\,\X}(i)\to\mathcal I_{T,\,C}(i)\to 0. 
 \end{gather*} Since $T$ is finite, $\h^1(\mathcal O_T(i))=0$.  So Serre's 
 computation of $\h^q(\mathcal O_{\X}(i))$ yields the formulas.  \end{proof} 

 \begin{theorem} Let $C\subset\X$ be a reduced curve of degree $d$. 
 Assume either the characteristic is $0$ or it is positive and does not 
 divide $d$.  Assume $C$ is a leaf of a foliation of $\X$ of degree 
 $m$. Let $\Sigma$ be the singular locus of $C$, and set 
 $\sigma:=\reg(\Sigma)$ and $\rho:=\sigma-d+2$. Then 
         $$d\le\begin{cases}m+1,         &\text{if }\rho\le0;\\ 
                             m+1+\rho,   &\text{if }\rho>0. 
         \end{cases}$$ 
  Furthermore, if $S$ is finite and $d\ge 2m+2$, then $d=m+1+\rho$. 
 \end{theorem} 

 \begin{proof} Suppose $\rho\le0$.  Then $\mathcal I_{\Sigma,\,\X}$ is 
 $(d-2)$-regular; in particular, $\h^1(\mathcal I_{\Sigma,\,\X}(d-3))=0$.  So 
 $\h^1(\mathcal I_{\Sigma,\,C}(d-3))=1$ by (2.4.2).  It now follows from (2.2.2) 
 that $\h^1(\Og^1_C)=1$.  Hence \cite{EK2}, Corollary 4.5, yields 
         $$\h^1(\mathcal O_C(m-1))=0.\eqno(2.5.1)$$ 

 Suppose $\rho>0$.  Then $\sigma>d-2$.  So (2.4.2) applies with
 $i:=\sigma-1$ and $\epsilon := 0$. 
 Now, $\h^1(\mathcal I_{\Sigma,\,\X}(\sigma-1))=0$ by definition of $\sigma$. 
 Hence 
         $$\h^1(\mathcal I_{\Sigma,\,C}(\sigma-1))=0.\eqno(2.5.2)$$ 

 By hypothesis, $C$ is a leaf of a foliation $\eta\:\Og^1_{\X}\to\mathcal 
 O_{\X}(m-1)$.  So $\eta$ induces a map $\vf\:\Og^1_C\to\mathcal O_C(m-1)$.  In 
 fact, $\eta|C=\vf\beta_C$ where $\beta_C\:\Og^1_X|C\to\Og^1_C$ is the 
 standard surjection.  Let $S$ be the singular locus of $\eta$.  By 
 (2.1.1) and (2.1.2), the image of $\eta$ is $\mathcal I_{S,\,\X}(m-1)$.
Hence the image of
 $\vf$ is $\mathcal I_{S\cap C,\, C}(m-1)$. 

 Since $C$ is reduced, $\mathcal O_C(m-1)$ is torsion free.  Hence, $\vf$ 
 factors through $\Og^1_C/\hbox{Torsion}$, which is equal to 
 $\mathcal I_{\Sigma,\,C}(d-3)$ by 
 (2.2.2).  Thus there is a map $\mathcal I_{\Sigma,\,C}(d-3) 
 \to\mathcal O_C(m-1)$.  It is injective since no component of $C$ lies in 
 $S$. Its image is $\mathcal I_{S\cap C,\, C}(m-1)$ by the preceding paragraph. 
 In other words, there is an exact sequence, 
  $$0\to\mathcal I_{\Sigma,\,C}(d-3)\to\mathcal O_C(m-1)\to\mathcal O_{S\cap C}(m-1)\to0.
         \eqno(2.5.3)$$ 
  Twist by $\rho$, and take cohomology.  Then use (2.5.2) to obtain 
         $$\h^1(\mathcal O_C(m+\rho-1))=0.\eqno(2.5.4)$$ 

 By duality, $\h^1(\mathcal O_C(i))=\h^0(\mathcal O_C(d-3-i))$ for any $i$.  Hence 
 $\h^1(\mathcal O_C(i))=0$ if and only if $i\geq d-2$.  Therefore, (2.5.1) and 
 (2.5.4) yield the first assertion. 

 Assume now that $d\ge 2m+2$.  Then $d\ge m+2$. So the first assertion 
 yields $d\le m+1+\rho$.  It remains to prove $d\ge m+1+\rho$. 

 The hypothesis on the characteristic implies that $\Sigma$ is the base 
 locus of the polar system; see the end of \S~2.2.  Hence $\Sigma$ is 
 contained in the finite intersection $T$ of two members of the system.  Now, 
 given any two curves $E,\,F\subset\X$ of degrees $e,\,f$, their 
 intersection $Z$, if finite, has regularity $e+f-1$, because the ideal 
 $\mathcal I_{Z,\,\X}$ has the Koszul resolution 
         $$0\to\mathcal O_{\X}(-E-F)\to\mathcal O_{\X}(-E)\oplus\mathcal O_{\X}(-F)\to 
         \mathcal I_{Z,\,\X}\to0.$$ 
  Therefore, $\reg(T)=2d-3$.  However, $T\supset\Sigma$ and $T$ is finite; 
 whence, $\sigma\le \reg(T)$.  Hence $\rho:=\sigma-d+2\le d-1$.  Thus
 if $m=0$, then $d\ge m+1+\rho$. 

 Assume $m>0$.  Then $d\ge2m+2\ge m+3$.  But $d\le m+1+\rho$ as observed 
 just above.  So $d\le m+\sigma-d+3\le \sigma$.  Therefore, we may take 
 $i:=\sigma-2$ in (2.4.2).  Now, $\h^1(\mathcal I_{\Sigma,\,\X}(\sigma-2))\neq0$ 
 by definition of $\sigma$.  Hence $\h^1(\mathcal 
 I_{\Sigma,\,C}(\sigma-2))\neq0$.  The exactness of (2.5.3) means $\mathcal 
 I_{\Sigma,\,C}(d-3)=\mathcal I_{S\cap C,\, C}(m-1)$.  Twisting by $\rho-1$ 
 and taking cohomology, we conclude 
         $$\h^1(\mathcal I_{S\cap C,\, C}(m+\rho-2))\neq0.\eqno(2.5.5)$$ 

 Proposition 2.3 implies $S$  is $i$-regular for all $i\ge 2m$.  So 
$\h^1(\mathcal I_{S,\,\X}(i-1))=0$ for all $i\ge 2m$.  Assume $S$ is finite.  Then 
 $\mathcal I_{S,\,\X}$ has finite colength in $\mathcal I_{S\cap C,\,\X}$.  Hence 
          $$\h^1(\mathcal I_{S\cap C,\,\X}(i-1))=0 \text{ for all }i\ge 2m. 
         \eqno(2.5.6)$$ 

 Set $j:=m+\rho-2$.  Then $j+1\ge d-2\ge 2m$.   
 Form the twisted exact sequence of ideals 
    $$0\to\mathcal I_{C,\,\X}(j)\to\mathcal I_{S\cap C,\,\X}(j)\to\mathcal I_{S\cap C,\,C}(j) 
         \to0.$$ 
  The first term is equal to $\mathcal O_{\X}(j-d)$.  Take cohomology.  Then 
 (2.5.6) and (2.5.5) yield $\h^2(\mathcal O_{\X}(j-d))\neq0$.  Hence $j-d\le 
 -3$, or $d\ge m+1+\rho$.  \end{proof} 

 \section{Bounds on the total Tjurina number} 

 \begin{lemma} 
   Let $C\subset\X$ be a reduced curve of degree $d\ge2$.  Let $\tau$ be 
 its total Tjurina number, and $u$ its number of 
 non-quasi-homogeneous singularities.  Let $\Sigma$ be its 
 singular locus, and set $\sigma:=\reg \Sigma$.  Assume the characteristic 
 is $0$.  Then 
         $$\sigma\le d-2 + (\tau-u)/(d-1).$$ 
 \end{lemma} 

 \begin{proof} Let $M$ be a general polar of $C$.  Then $M$ is smooth off 
 the base locus of the polar system by Bertini's  First Theorem, Theorem 
 (3.2) in \cite{K2}, because the characteristic is $0$.  This base locus 
 is $\Sigma$, again because the characteristic is $0$; see the end of 
 \S~2.2.  But $\Sigma$ is finite, since $C$ is reduced.  Hence, $M$ is 
 generically smooth.  So, since $M$ is Cohen--Macaulay, $M$ is reduced. 

 Either $M$ is irreducible or $C$ is the union of concurrent lines by 
 \cite{GL}, Lemma 3.8, p.~333.  An alternative proof of this fact is 
 given at the end of the present proof. 

 First, assume $C$ is the union of concurrent lines, say meeting at $P$. 
 Then $P$ is a quasi-homogeneous singularity of $C$, and its only 
 singularity; so $u=0$.  Since $C$ is a cone, its polar system has 
 dimension 1.  Hence the base locus $\Sigma$ is the complete intersection 
 of two polars.  So $\tau=(d-1)^2$ by Bezout's theorem.  In addition, 
 $\sigma=2d-3$; see the middle of the proof of Theorem 2.5.  Thus, in the 
 present case, equality holds in the asserted inequality. 

 {}From now on, assume $M$ is irreducible. The lemma holds if $\sigma\leq 
 d-2$ since $\tau-u\geq 0$; so we may assume $\sigma\geq d-1$. 

 Since $\sigma\geq1$, we have $\h^2(\mathcal I_{\Sigma,\,\X}(\sigma-3))=0$ by 
 (2.4.1).  Hence $\h^1(\mathcal I_{\Sigma,\,\X}(\sigma-2))\neq0$ by definition of 
 $\sigma$.  So (2.4.2) yields $\h^1(\mathcal I_{\Sigma,M}(\sigma-2))\neq0$. 
 Hence duality yields 
         $$\Hom\bigl(\mathcal I_{\Sigma,M}(\sigma-2),\,\mathcal O_M(d-4)\bigr)\neq0. 
         \eqno(3.1.1)$$ 

 Given $P\in \Sigma$, denote by $\ve_P$ the multiplicity of the ideal
of $\Sigma$  in the local ring $\mathcal O_{\X,\,P}$.  Then two general
polars $M$ and $N$ meet at $P$ with multiplicity $\ve_P$.  On the one
hand, the Tjurina number $\tau_P$ of $C$ at $P$ is the length at $P$ of
the scheme cut out by all the polars.  On the other hand, the Milnor
number $\mu_P$ of $C$ at $P$ is equal to the intersection number at $P$
of two special polars.  Hence
	$$\tau_P\le\ve_P\le\mu_P.\eqno(3.1.2)$$ 
  Saito's theorem (see \S~2.2) asserts  $\tau_P=\mu_P$ if and only if $P$ 
 is quasi-homogeneous. 

 Suppose $\mathcal I_{\Sigma,M}$ is invertible at $P\in \Sigma$.  Let's prove 
 $P$ is quasi-homo\-ge\-ne\-ous on $C$.  (We don't need the converse, but 
 it holds as $M$ is general.  Indeed, if $\tau_P=\mu_P$, then 
 $\tau_P=\ve_P$ by (3.1.2), and the equation $\tau_P=\ve_P$ just means 
 $\Sigma$ is cut out of $\X$ at $P$ by two general polars.) 

 Take local coordinates $x,\,y$ for $\X$ at $P$, and let $f(x,y)=0$ be a 
 local equation for $C$.  Then $\mathcal I_{\Sigma,\,\X}$ is generated at $P$ by 
 $f$ and its partials $f_x,\,f_y$.  If the latter alone generate, then 
 $\tau_P=\mu_P$, and so $P$ is quasi-homo\-ge\-ne\-ous on $C$ by Saito's 
 theorem, recalled above.  Suppose $f_x,\,f_y$ don't generate, and let's 
 achieve a contradiction. 

 At any rate, $\mathcal I_{\Sigma,\X}$ is 2-generated at $P$ as $\mathcal 
 I_{\Sigma,M}$ is invertible there.  Hence $\mathcal I_{\Sigma,\,\X}$ is 
 generated at $P$ either by $f,\,f_x$ or by $f,\,f_y$.  Either way, $\mathcal 
 I_{\Sigma,\,C}$ is invertible at $P$.  But then, (2.2.2) implies 
 $\Hom(\Og^1_C,\mathcal O_C)$ is invertible at $P$.  So $P\notin \Sigma$ by 
 \cite{L}, Theorem~1, p.~879, a contradiction. 

 Alternatively, since $\mathcal I_{\Sigma,\,C}$ is invertible at $P$, it follows 
 that $\tau_P=e_P$ where $e_P$ denotes the multiplicity of $\mathcal I_{\Sigma,\,C}$ 
 at $P$.  However, Teissier proved $e_P=\mu_P+m_P-1$ where $m_P$ denotes the 
 multiplicity of $C$ at $P$; see Corollaire 1.5 on p.~320 and Remarque 
 1.6~1) on p.~300 in \cite{Te}; alternatively, see pp.~358--359 of 
 \cite{K1}, where Teissier's formula is derived from the Milnor--Jung 
 formula and another formula for $e_P$, due to Piene.  Since 
 $\tau_P\le\mu_P$ by (3.1.2), it follows that $m_P=1$.  So again $P\notin 
 \Sigma$, a contradiction. 

 Owing to (3.1.1), there is a nonzero map $w\:\mathcal I_{\Sigma,M}\to\mathcal 
 O_M(d-\sigma-2)$.  Since $M$ is reduced and irreducible, $w$ is 
 injective and its cokernel $\Cok w$ is supported by a finite subset 
 $W\subset M$.  If $P\in C$ is one of the $u$ singularities that are not 
 quasi-homogeneous, then $\mathcal I_{\Sigma,M}$ is not invertible at $P$, as 
 we just proved.  Hence $w$ is not surjective at $P$; in other words, 
 $P\in W$.  It follows that 
  \begin{align*} 
  u&\le\deg W\le\deg\Cok w=\chi(\mathcal O_M(d-\sigma-2))-\chi(\mathcal I_{\Sigma,M})\\ 
   &=\bigl(\chi(\mathcal O_M)-\chi(\mathcal I_{\Sigma,M})\bigr) 
    -\bigl(\chi(\mathcal O_M)-\chi(\mathcal O_M(d-\sigma-2))\bigr)\\ 
   &=\tau-(d-1)(\sigma-d+2). 
  \end{align*} 
  The asserted bound on $\sigma$ follows directly, and the proof is 
 complete. 

 Here's an alternative proof that either $M$ is irreducible or $C$ is the 
 union of concurrent lines.  The Milnor--Jung formula says 
 $\mu_P=2\delta_P-(r_P-1)$ where $\delta_P$ is the genus diminution at 
 $P$ and $r_P$ is the number of branches of $C$ at $P$.  Set 
 $\delta:=\sum\delta_P$, and let $r$ be the number of irreducible 
 components of $C$. Since $C$ is connected, $\sum(r_P-1)\geq 
 r-1$. Hence, 
    	$$\sum\mu_P\leq 2\delta-(r-1). \eqno(3.1.3)$$ 

 Let $p_a$ be the arithmetic genus of $C$, and $g$ its geometric genus. 
 Then $\delta=p_a-g+r-1$. In addition, $p_a=(d-1)(d-2)/2$ because 
 $C\subseteq\IP^2$ is of degree $d$. Since $g\geq 0$, we get 
  	$$2\delta-(r-1)\leq (d-1)(d-2)+(r-1).$$ 
  But $r\leq d$.  So the above inequality and (3.1.3) yield 
  	$$\sum\mu_P\leq (d-1)^2.  \eqno(3.1.4)$$ 

  By Bezout's theorem, two general polars $M$ and $N$ meet in $(d-1)^2$ 
 points  $P$  counted with multiplicity.  If $P\in\Sigma$, then this 
multiplicity is  $\ve_P$.  Hence 
  	$$\textstyle \sum_{P\in\Sigma}\ve_P\leq (d-1)^2.  \eqno(3.1.5)$$ 
  In addition, either equality holds in (3.1.5) or $M$ and $N$ meet outside 
of  $\Sigma$. 

 First, assume equality holds in (3.1.5).  Then equality holds in 
(3.1.4) because of (3.1.2). Hence also $r=d$ and 
 $\sum(r_P-1)=r-1$.  The first equation says $C$ is a union of 
 lines.  Then the second says these lines are concurrent. 

 Finally, assume $M$ and $N$ meet outside $\Sigma$.  Then the polar 
 system is not composite with a pencil.  Hence, by Bertini's second 
 theorem, Theorem (5.1) in \cite{K2}, a general polar is irreducible. 
   \end{proof} 

  \begin{theorem}[du Plessis--Wall] Let $C\subset\X$ be a reduced curve of 
 degree $d\ge2$.  Let $\tau$ be its total Tjurina number, and $u$ its 
 number of non-quasi-homogeneous singularities.  Let $\Sigma$ be its 
 singular locus, and set $\sigma:=\reg \Sigma$.  Assume the 
 characteristic is $0$.  Assume $C$ is a leaf of a foliation of $\IP^2$ 
 of degree $m$.  Then 
         $$(d-1)(d-m-1)+u\le\tau.$$ 
 If equality holds, then either $d=m+1$ and $C$ is smooth, or $d>m+1$ 
 and $\sigma=2d-m-3$. 
 \end{theorem} 

 \begin{proof} 
 First, assume $d\leq m+1$. Then the asserted inequality holds because 
 $\tau\geq u$. Furthermore, equality holds if and only if $d=m+1$ and 
 $\tau=u$. Now, $\tau=u$ if and only if every singularity of $C$ is 
 non-quasi-homogeneous and has Tjurina number 1.  However, only an 
 ordinary node has Tjurina number 1, and it is quasi-homogeneous.  Thus 
 $\tau=u$ if and only if $\tau=0$. 

 Finally, assume $d>m+1$.  Then $d-m-1\leq\sigma+2-d$ by Theorem 2.5.  In 
 addition, $(d-1)(\sigma+2-d)\leq\tau-u$ by Lemma 3.1. The asserted 
 inequality follows directly.  Furthermore, if equality holds, then 
 $d-m-1=\sigma+2-d$. 
  \end{proof} 

 \begin{proposition}[du Plessis--Wall] Let $C\subset\X$ be a reduced curve of 
 degree $d$ and total Tjurina number $\tau$.  Let $m$ be the least degree 
 of a foliation of\/ $\X$ with $C$ as leaf.  Then 
         $$m\le d-1 \hbox{ and }\tau\le(d-1)(d-m-1)+m^2.$$ 
 Moreover, if $d\le 2m$ and if $C$ is irreducible, then 
         $$\textstyle\tau\le(d-1)(d-m-1)+m^2-\binom{2m+2-d}{2}.$$ 
 \end{proposition} 

  \begin{proof} To see $m\le d-1$, take homogeneous coordinates 
 $x,\,y,\,z$ for $\X$, and let $F(x,y,z)=0$ define $C$. 
 Given a component of $C$, there is a partial of $F$ that doesn't vanish 
 along it since $C$ is reduced.  So some linear combination of the 
 partials doesn't vanish along any component.  Changing coordinates, we 
 may assume neither $F_x$ nor $z$ vanishes along any component. 

Form the Hamilton foliation 
$\theta\:\Og^1_{\X}\to\mathcal O_{\X}(-1)^3\to\mathcal O_{\X}(d-2)$ where the second 
map is $(F_y,-F_x,0)$.  Then $\theta|C$ 
 factors through $\beta_C\:\Og^1_{\X}|C\to\Og^1_C$ because 
 $F_yF_x=F_xF_y$.  Suppose $\theta$ vanishes along some component $C'$ of 
 $C$.  Then $(F_y,-F_x,0)$ is a polynomial multiple of $(x,y,z)$ on $C'$ 
 since the restriction of the Euler sequence $0\to\Og^1_{\X}\to\mathcal 
 O_{\X}(-1)^3\to\mathcal O_{\X}\to0$ is exact where the third map is $(x,y,z)$. 
 Hence, either $z$ or $F_x$ must vanish along $C'$ , but neither do. 
 Thus $C$ is a leaf of $\theta$.  Now, $\theta$ is of degree $d-1$. 
 Therefore $m\le d-1$.

 The first bound on $\tau$ is proved for curves in higher space in 
 \cite{EK1}, Proposition 6.3.  In the present case, the proof becomes 
 shorter.  We give it now mainly to recall its ingredients. 

 Let $\eta\:\Og^1_{\X}\to\mathcal O_{\X}(m-1)$ be a foliation with $C$ as leaf, 
 and $S$ its singular locus. Since $m$ is minimal, $S$ is finite. Let 
 $\Sigma$ be the singular locus of $C$. The exactness of (2.5.3) means 
 $\mathcal I_{\Sigma,\,C}(d-3)=\mathcal I_{S\cap C,\, C}(m-1)$.  So $\mathcal 
 I_{\Sigma,\,C}(d-m-2)=\mathcal I_{S\cap C,\, C}$.  Hence 
  \begin{align*} 
   \tau&=\chi(\mathcal O_C)-\chi(\mathcal I_{\Sigma,\,C}) 
         =\chi(\mathcal O_C(d-m-2))-\chi(\mathcal I_{\Sigma,\,C}(d-m-2))\\ 
         &=\bigl(\chi(\mathcal O_C(d-m-2)-\chi(\mathcal O_C)\bigr) 
           +\bigl(\chi(\mathcal O_C)-\chi(\mathcal I_{S\cap C,\, C})\bigr)\\ 
         &=d(d-m-2)+\deg(S\cap C)\tag{3.3.1} 
  \end{align*} 

 Since $S$ is finite, $S$ represents the top Chern class of
 $(\Og^1_{\X})^*(m-1)$; so 
         $$\deg S=m^2+m+1.\eqno(3.3.2)$$
 Alternatively, $\deg S$ can be computed from Sequence (2.3.1).  The
first assertion now follows from (3.3.1), from (3.3.2) and from the
trivial bound $\deg(S\cap C)\le\deg S$.

 Better lower bounds for $\deg S-\deg(S\cap C)$ yield better upper bounds 
 for $\tau$.  For example, if $d\le 2m$ and if $C$ is irreducible, then 
 the following bound obtains: 
   $$\textstyle\deg S-\deg(S\cap C)\ge\binom{2m+2-d}{2}.\eqno(3.3.3)$$ 
  And it,  (3.3.1), and (3.3.2) yield the second assertion. 
 This bound is proved next. 

 Assume $d\le 2m$.  Then $\mathcal I_{C,\,\X}=\mathcal O_{\X}(-d)$ yields $\h^0(\mathcal 
I_{C,\,\X}(2m))= \binom{2m+2-d}{2}$.  Now, $\mathcal I_{C,\,\X}$ generates 
$\mathcal I_{S\cap 
 C,\,S}$ in $\mathcal O_S$.  Hence (3.3.3) obtains if the following 
 composition is injective: 
 $$ 
u\:\H^0(\mathcal I_{C,\,\X}(2m))\to\H^0(\mathcal O_{\X}(2m))\to\H^0(\mathcal O_S(2m)). 
 $$ 

 Extend the exact sequence (2.3.1) to an augmented resolution: 
  $$0\to\Og^2_{\X}(-2m+2)\to\Og^1_{\X}(-m+1)\to\mathcal O_{\X}\to\mathcal O_S\to 0.$$ 
 Twist it by $2m$ and take global sections to obtain this complex: 
         $$\H^0(\Og^1_{\X}(m+1))\to\H^0(\mathcal O_{\X}(2m))\to\H^0(\mathcal O_S(2m)). 
         \eqno(3.3.4)$$ 
 It is exact because $\H^1(\Og^2_{\X}(2))=0$. 

 Given $s\in\H^0(\mathcal I_{C,\,\X}(2m))$ with $u(s)=0$, we must show $s=0$. 

  Since (3.3.4) is exact, $s$ is the image of some $s'\in 
 \H^0(\Og^1_{\X}(m+1))$.  Let $\pi\:\Og^1_{\X}\ox\Og^1_{\X}\to\Og^2_{\X}$ 
 be the natural pairing, and consider the composition, 
 $$ 
 \begin{CD} 
 \eta'\:\Og^1_{\X} @>1\ox s'>> \Og^1_{\X}\ox\Og^1_{\X}(m+1) @>\pi(m+1)>> 
 \Og^2_{\X}(m+1). 
 \end{CD} 
 $$ 
 It turns out that $\eta'$ is a foliation with $C$ as leaf if $s\neq 0$ 
 and if $C$ is irreducible.  But then $\eta'$ has degree $m-1$, 
 contradicting the minimality of $m$.  Thus necessarily $s=0$, as 
 desired. 

 So assume $s\neq 0$.  Then $s'\neq 0$.  Also $\eta'\neq 0$ 
 because $\pi$ is a perfect pairing.  Thus $\eta'$ is a foliation.  And it 
 remains to show $C$ is a leaf of $\eta'$. 

 Let $\mathcal K$ be the kernel of $\beta_C\:\Og^1_{\X}|C\to\Og^1_C$.  Since 
 $C$ is reduced, $C$ is generically smooth and $\mathcal O_C(-d)$ is torsion 
 free.  So the first map in (2.2.1) is injective.  Thus $\mathcal K=\mathcal 
 O_C(-d)$. 

 Note $s'|C\in\H^0(\Ker(\eta(m+1)|C))$ because $s\in\H^0(\mathcal I_{C,\,\X}(2m))$. 
 Since $C$ is a leaf of $\eta$, there is, by definition, a map 
 $\vf\:\Og^1_C\to\mathcal O_C(m-1)$ such that $\eta|C=\vf\beta_C$.  Since $S$ is 
 finite and $C$ is generically smooth, $\vf$ is generically an 
 isomorphism. Thus $\Ker(\eta(m+1)|C)$ and $\mathcal K(m+1)$ agree generically, 
 and hence $\Im(s'|C)$ generically lies in $\mathcal K(m+1)$. 

 It follows that $\eta'|C$ factors through $\beta_C$, or 
 $(\eta'|C)\mathcal K=0$.  Indeed, since $\Im(s'|C)$  generically 
 lies in $\mathcal K(m+1)$, also $(\eta'|C)\mathcal K$ generically lies in 
 $(\pi|C)(m+1)(\mathcal K\ox\mathcal K(m+1))$.  But $\mathcal K$ is invertible; whence, 
 $(\pi|C)(\mathcal K\ox\mathcal K)=0$. So $(\eta'|C)\mathcal K$ is generically trivial. 
But $(\eta'|C)\mathcal K$ lies in $\Og^2_{\X}|C(m+1)$, and $\Og^2_{\X}|C(m+1)$ 
is torsion 
 free.  Therefore, $(\eta'|C)\mathcal K=0$. 

 Finally, $C$ does not lie in the singular locus of $\eta'$; otherwise, 
$\eta'$ would factor through a foliation of degree $m-1-d$, but 
$m-1-d<0$.  So, since $C$ is irreducible, no component of $C$ lies in 
 the singular locus.  Thus $C$ is indeed a leaf of $\eta'$. 
 \end{proof} 

  \begin{remark} Over $\bf C$, let $C\subset\X $ be a reduced curve of
degree $d$ and with total Tjurina number $\tau$.  Let $m$ be the least
degree of a foliation $\eta$ of $\X $ such that $\eta|C$ factors through
the natural surjection $\beta_C\:\Og^1_{\X }|C\to\Og^1_C$.  In
essentially this setting, Du Plessis and
 Wall, in \cite{dPW}, Theorem 3.2, p.~263, proved $d\ge m+1$ and 
	$$(d-1)(d-m-1)\le\tau\le(d-1)(d-m-1)+m^2;\eqno(3.4.1)$$ 
 furthermore, if $d\le 2m$, then 
	$$\textstyle\tau\le(d-1)(d-m-1)+m^2-\binom{2m+2-d}{2}.\eqno(3.4.2)$$ 

Let $m'$ be the least degree of a foliation of $\X$ with $C$ as leaf. 
Clearly, $m\le m'$. If $m=m'$, then (3.4.1) follows 
from Theorem 3.2 and Proposition~3.3. In this case, Theorem 3.2 is stronger, 
since $u\ge 0$. If also $C$ is irreducible, then (3.4.2) 
follows from Proposition~3.3. 

In fact, if $C$ is irreducible, then $m=m'$. Indeed, let $\eta$ be a 
foliation of degree $m$ of $\X$ such that $\eta|C$ factors through 
$\beta_C$. If $C$ were not a leaf of $\eta$, then $C$ would be contained in 
the singular locus of $\eta$, and hence $m\ge d$. But $m\le d-1$ by 
Proposition~3.3.

However, if $C$ is reducible then possibly $m<m'$.
 For instance, suppose $C$ is the union of $d$ lines, $d-1$ concurrent. 
Let $K$ be the union of the latter, $M$ the additional line.  
Since $C$ is not a cone, $m>0$.  But $K$ is a leaf of a foliation of 
degree 0.  So there is a foliation $\eta$ of degree 1, with $M$ in its 
singular locus, and such that $\eta|C$ factors through $\beta_C$.  
Thus $m=1$.  However, given a foliation 
$\eta'$ of $\X$ of degree $m'$ with $C$ as leaf, it follows from
\cite{dPW}, Corollary 3.1.1, p. 263, that either $m+m'\ge d-1$ or both 
$\eta$ and $\eta'$ factor through a foliation $\eta''$ such that 
$\eta''|C$ factors through $\beta_C$. If the latter were true, then 
$\eta$ would be a scalar multiple of $\eta''$, because $m$ is minimal. 
Then $\eta'$ would factor through $\eta$, and hence the singular locus 
of $\eta$ would be contained in that of $\eta'$. But $M$ is in the singular 
locus of $\eta$ and not in that of $\eta'$. So $m+m'\ge d-1$, and hence 
$m'\ge d-2$.
 
 It is possible to recover Du Plessis and Wall's lower bound on $\tau$ 
 with a bit more work.  Given a foliation $\eta$ such that $\eta|C$ 
 factors through $\beta_C$, let $B\subseteq C$ be the union of all the 
 1-dimensional components of the singular locus of $\eta$ contained in 
 $C$.  Then $\eta$ factors through a foliation $\eta'$ whose degree is 
 $d-\deg B$ and for which $A:=C-B$ is a leaf.  Then Theorem 3.2 applies 
 to $A$ and $\eta'$ in place of $C$ and $\eta$.  The resulting lower 
 bound on the Tjurina number of $A$ can now be used to obtain the desired 
 lower bound on $\tau$.  This argument is described in detail in the 
 proof of \cite{EK1}, Corollary~6.4.  However, the argument does not seem 
 to yield the upper bounds. 
 \end{remark}

 \end{document}